\documentclass{article}

\usepackage{arxiv}

\usepackage[utf8]{inputenc} 
\usepackage[T1]{fontenc}    
\usepackage{hyperref}       
\usepackage{url}            
\usepackage{booktabs}       
\usepackage{amsfonts}       
\usepackage{nicefrac}       
\usepackage{microtype}      
\usepackage{lipsum}
\usepackage{amsmath,amsthm,amssymb}
\usepackage[all]{xy}
\usepackage{graphicx}
\usepackage{mathrsfs}
\usepackage{ upgreek}
\theoremstyle{plain}
\theoremstyle{remark}

\newtheorem*{convention*}{Convention}
\theoremstyle{plain}
\newtheorem{theorem}{Theorem}[section]

\theoremstyle{definition}

\newtheorem{example}[theorem]{Example}

\title{A  Mathematical  Model for Arch Fingerprint}

\author{Ibrahim Jawarneh\\
  Department of Mathematics\\
 Al-Hussein Bin Talal University\\
Ma'an, P.O. Box (20), 71111, Jordan\\
  \texttt{ibrahim.a.jawarneh@ahu.edu.jo} \\
   \And
  Nesreen Alsharman\\
  Computer Science\\
The World Islamic Sciences and Education University\\
 Amman, Jordan\\
  \texttt{nesreen.alsharman@wise.edu.jo} \\
}

\begin{document}
\maketitle

\begin{abstract}
In this paper, different categories of the arch fingerprint are modelled in a general dynamical system  using ordinary differential equations with a parameter $\theta > 0$. We study its global dynamics and analyze the existence and stability of equilibria. Numerical simulations using Maple show the matching between real images of categories of arch fingerprint and phase portraits of the considered dynamical system.
\end{abstract}


\section{Introduction}\label{FP introduction}
Fingerprints are patterns, made by friction ridges of a human finger, which appear on the pads of the fingers and thumbs. The fingerprints is the most accurate and reliable for identifying person, since fingerprints are most unique biometric characteristics to any person; therefore it is used in forensic divisions worldwide for security, and in criminal cases where even the twins have non-identical fingerprints. Early study  about fingerprints appeared in 1892 by Sir Francis Galton in his book, finger prints \cite{Galton1892}. The  patterns of fingerprints are classified into three general types; loop, whorl, and arch, many studies indicate that arch patterns are observed in about 5\% of all fingerprints, see \cite{Arent2019,Galton1892,Henry1990}. 
\par 
\indent The arch  patterns  display a relatively horizontal ridges run from the left to the right side of the fingerprint with growth in the center.  This paper discusses three categories of arch fingerprint which are shown in figure \ref{Plain, Tented and strong arch fingerprint} and we depend on it as a database of the classes of arch fingerprint.
\begin{itemize}
\item  Plain arch in which the ridges flow relatively horizontally with a little rise in the middle where delta is not clear here,  see   picture (a) in figure \ref{Plain, Tented and strong arch fingerprint}.
\item Tented arch is similar to the plain arch   with a bigger rise, and  at least one ridge  with  short length is vertically oriented in the middle where delta is clear here, see picture (b) in figure \ref{Plain, Tented and strong arch fingerprint}.
\item Strong arch is similar to the tented arch with many relatively longer ridges are vertically oriented in the middle where delta is very clear here, this category looks like Christmas tree, see picture (c) in  figure \ref{Plain, Tented and strong arch fingerprint}.
\end{itemize}
\begin{figure}[ht]
\centering
\begin{minipage}[c]{0.33\linewidth}
\centering
\includegraphics[width=1.9in]{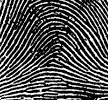}\\ \small (a)
\end{minipage}%
\begin{minipage}[c]{0.33\linewidth}
\centering
\includegraphics[width=2.1in]{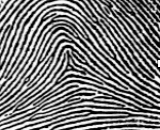}\\ \small (b)
\end{minipage}%
\centering
\begin{minipage}[c]{0.33\linewidth}
\centering
\includegraphics[width=1.9in]{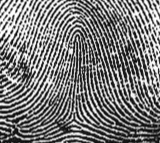}\\ \small (c)
\end{minipage}%
\caption{(a) Plain arch fingerprint,  (b) tented arch fingerprint and (c) strong arch fingerprint}
\label{Plain, Tented and strong arch fingerprint}
\end{figure}
\par
\indent A few studies talked about modelling of arch fingerprint, Huckemann et al. \cite{Huckemann2008} proposed a model based on quadratic differentials that describes arch fingerprint but they  do not show the results for a tented arch since they considered it as a special case of a loop. Moerover, they did not study the strong arch.  The conditional cosine functions were used for modeling arch-typed  fingerprints in the full plane without going over specific classes, see \cite{Cappelli2007, Maltoni2009,Wang2011}. We did not find  any study  proposed a general mathematical model for all classes of arch fingerprint as we focus in this paper. Fingerprints can be captured as graphical ridge and valley patterns, so we start to put a model which simulates the  ridge flow and pattern types in the above classes. Using differential equations in creating dynamical system that describes plain, tented, and strong arch fingerprint require understanding the behaviour of the ridges in these classes and how much the delta is clear in each class. In this case the flow of the phase portrait of this system represents the ridges in the fingerprint, so we think about the deformation of the phase portrait of straight flow to generate such system. We proceeded by trial and error to approximate the following dynamical system that represents the general map of the above categories of the arch. 
\begin{equation}\label{FP-Arch-sys}
 \begin{aligned} 
  \dot{x} &= y^2,
\\
\dot {y}  &= -{\theta}x, \quad  \theta > 0.
\end{aligned}   
\end{equation}
\par
This paper is organized as follows :  In the next section, we study the stability of the equilibria of system \eqref{FP-Arch-sys}.  In section \ref{FP-Arch Discussion} we go over different values of $\theta$ in the dynamical system \eqref{FP-Arch-sys} to simulate the  categories of the arch fingerprint. In section \ref{FP-Arch conculusion} we conclude the results of the proposed  model.
\section{Steady States and  Their Stability} \label{FP-Arch staeady state}
\hspace{\parindent}
To study the stability of the system $\eqref{FP-Arch-sys}$, we find the equilibria first, which are the solutions of the following equations:
\begin{align} 
    0  & = y^2,  \label{Arch FP 1st order 1st eqn}
    \\
    0  & = -  {\theta}x , \label{Arch FP 1st order 2nd eqn}
\end{align}  
and are given by $E_0 =(0,0)$ which is called equilibrium point or singular point. The Jacobian matrix of $\eqref{FP-Arch-sys}$ takes the form
\begin{equation} \label{Jacobian Arch}
J =
  \begin{bmatrix}
  0 & 2y \\
-\theta & 0
\end{bmatrix}  .
\end{equation}
The Jacobian matrix evaluated at the equilibrium point $E_0 = (0,0)$ is
\begin{equation} \label{JacobianE_0 Arch}
J(E_0) =
  \begin{bmatrix}
  0 &  0\\
-\theta & 0
\end{bmatrix}  .
\end{equation}
From the Jacobian matrix \eqref{JacobianE_0 Arch}, the eigenvalues of $J(E_0)$ are $\lambda_{1,2} = 0$ which means degenerate nonhyperbolic equilibrium point,
the phase portrait for this system is shown in section \ref{FP-Arch Discussion}. In all categories of arch fingerprints, we see that a deleted neighborhood of the origin consists of upper and lower hyperbolic sectors, also it consists of left and right separatrices. This type of critical point is called a cusp, see \cite{Perko2001}. The cusp  plays the role of the delta  in the above description of the categories of arch fingerprint. But the flow which is located above the origin makes different types of the angle depending on the value of $\theta$, for example we have obtuse angle in plain category, around right angle  in tented category, and acute angle in strong category.  
\section{Simulations and Numerical Results} \label{FP-Arch Discussion}
\hspace{\parindent}
In this section, we study the system \eqref{FP-Arch-sys} with different values of $\theta$ and display its simulations and their matching images of the above categories of the arch fingerprint. We use bold red lines to represent the separatrices, the green lines at different initial conditions above x-axis to show the above hyperbolic sector, and the brown line at different initial conditions below x-axis to display the below  hyperbolic sector.
\subsection{Plain arch fingerprint }
\begin{figure}[ht]
\centering
\begin{minipage}[c]{0.38\linewidth}
\centering
\includegraphics[width=2.3in]{plain3.png}
\end{minipage}%
\begin{minipage}[c]{0.38\linewidth}
\centering
\includegraphics[width=2.2 in]{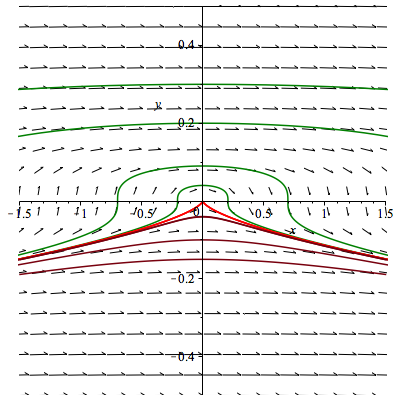}
\end{minipage}%
\caption{Simulation of the plain fingerprint}
\label{Plain fingerprint}
\end{figure}

\begin{figure}[ht]
\centering
\begin{minipage}[c]{0.4\linewidth}
\centering
\includegraphics[width= 2.3 in]{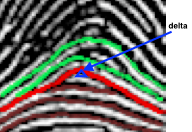}\\ \small (a)
\end{minipage}%
\begin{minipage}[c]{0.4\linewidth}
\centering
\includegraphics[width= 2.5 in]{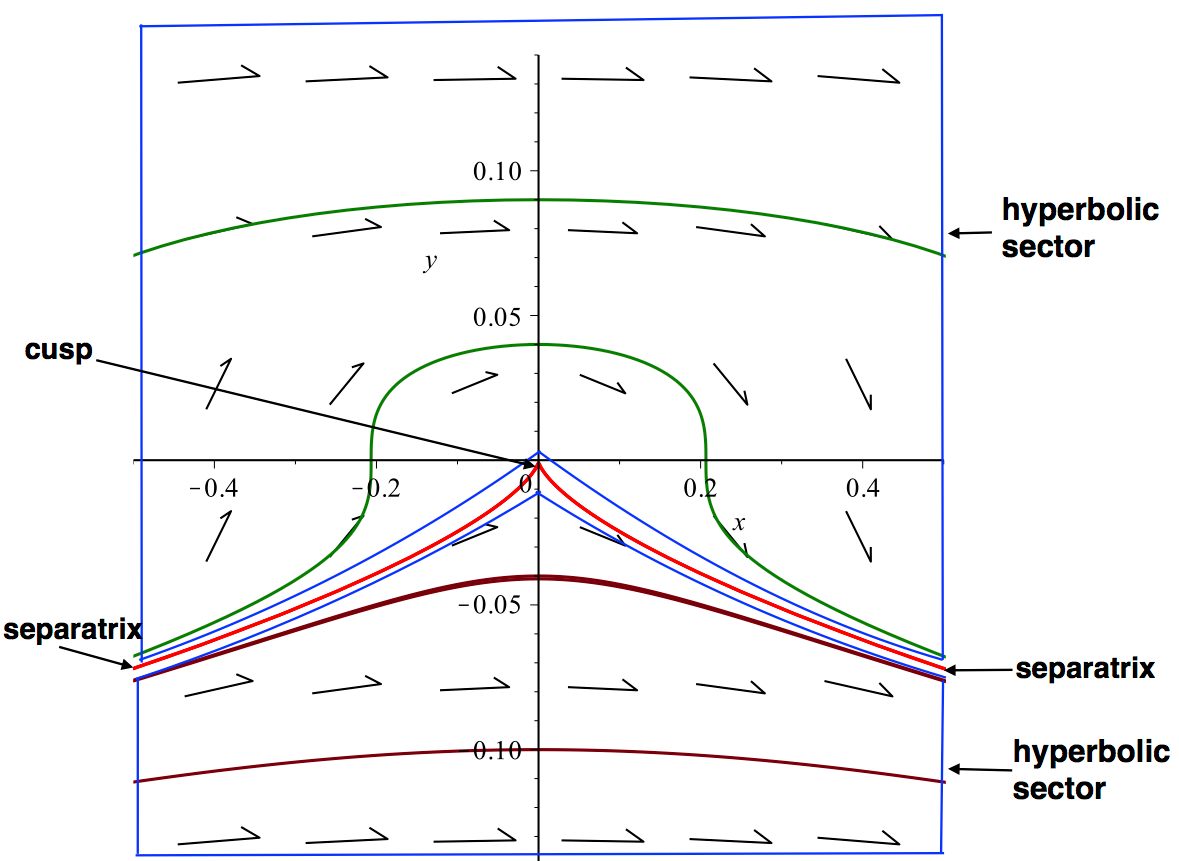}\\ \small (b)
\end{minipage}%
\caption{(a) Neighborhood of the delta in the image of the plain arch fingerprint and (b) Neighborhood of the origin in the phase portrait of example  \ref{example Plain arch}   .}
\label{Nhd Plain}
\end{figure}

\hspace{\parindent} When $\theta$ is very small and closed to zero in \eqref{FP-Arch-sys}, the flow of phase portrait is relatively horizontally with a slight growth in the middle as in figure \ref{Plain fingerprint}. For example consider system \eqref{FP-Arch-sys} with  $\theta = 0.001$.
\begin{example} \label{example Plain arch}
\begin{equation} \label{theta = 0.001}
   \begin{aligned} 
      \dot{x} &= y^2,
\\
\dot{y} &= -0.001x.
\end{aligned} 
\end{equation}
\end{example}
Figure \ref{Plain fingerprint} shows the simulation of the example \eqref{example Plain arch}  using Maple.
Figure \ref{Nhd Plain} explains a closer look at the neighborhood of the origin, we have explained the separatrices in the phase portrait, and we  use a thin blue line to determine the region  of hyperbolic sectors, and hence the origin is cusp fixed point which represents the delta in the plain arch fingerprint image. We notice the similarity between the flow in the phase portrait of the system \eqref{theta = 0.001} and the layout of the ridges  in the image of the plain arch fingerprint.

\subsection{Tented arch fingerprint }
\begin{figure}[ht]
\centering
\begin{minipage}[c]{0.38\linewidth}
\centering
\includegraphics[width=2.3 in]{tented.png}
\end{minipage}%
\begin{minipage}[c]{0.38\linewidth}
\centering
\includegraphics[width=2.1 in]{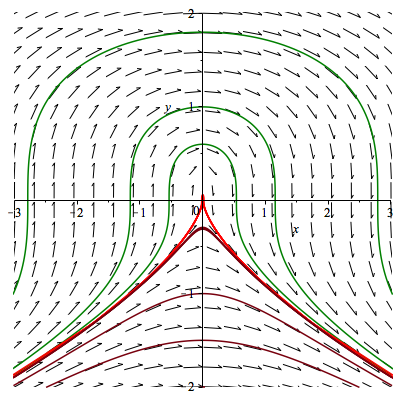}
\end{minipage}%
\caption{ Simulation of the tented arch fingerprint}
\label{Tented fingerprint}
\end{figure}

\begin{figure}[ht]
\centering
\begin{minipage}[c]{0.4\linewidth}
\centering
\includegraphics[width= 2.3 in]{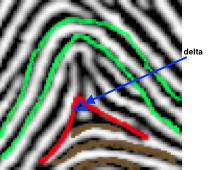}\\ \small (a)
\end{minipage}%
\begin{minipage}[c]{0.4\linewidth}
\centering
\includegraphics[width= 2.6 in]{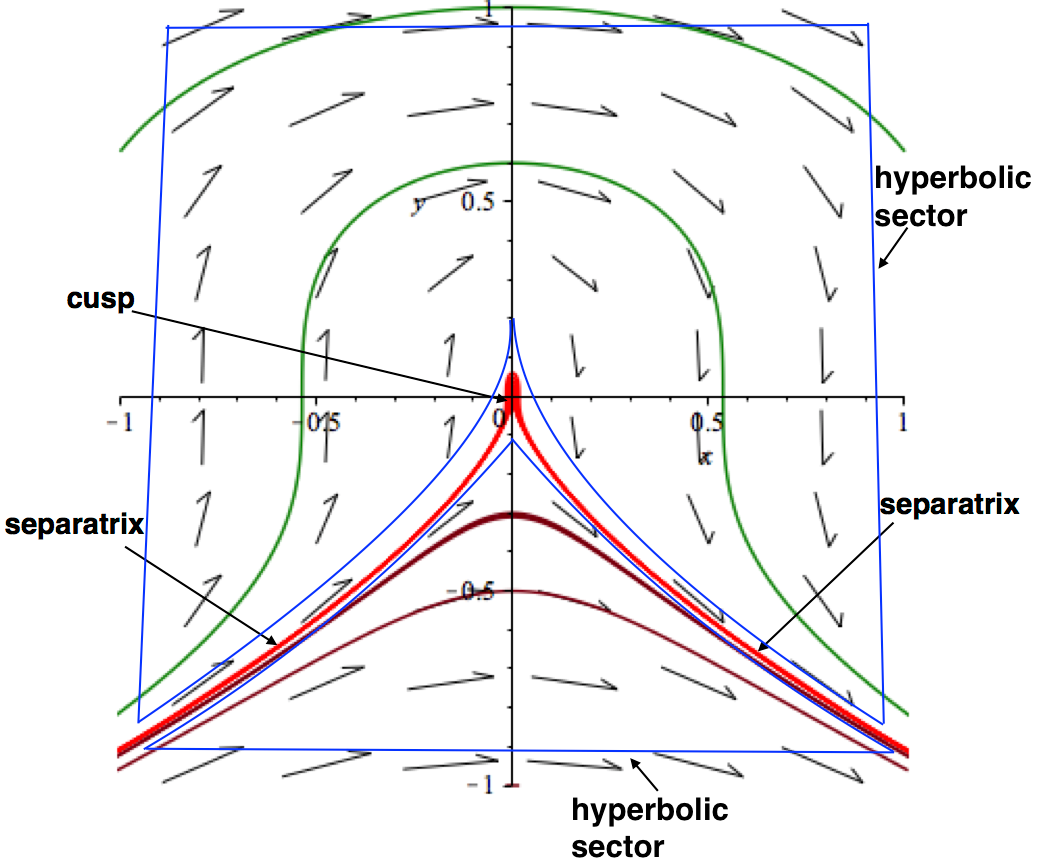}\\ \small (b)
\end{minipage}%
\caption{(a) Neighborhood of the delta in the image of the tented arch fingerprint and (b) Neighborhood of the origin in the phase portrait of example  \ref{example tented arch}   .}
\label{Tented NBH fingerprint}
\end{figure}

\hspace{\parindent}
If we increase  $\theta$   reaching around  $0.5$ in \eqref{FP-Arch-sys},  the flow rises more in the middle with a smaller angle  above the origin than the case of plain type. In this case we get phase portrait looks like tented arch, for example consider system \eqref{FP-Arch-sys} with  $\theta = 0.5$.
\begin{example} \label{example tented arch}
\begin{equation} \label{theta = 1}
   \begin{aligned} 
      \frac{dx}{dt} &= y^2,
\\
\frac{dy}{dt}  &= - 0.5x.
\end{aligned} 
\end{equation}
\end{example}
The phase portrait of the example \eqref{example tented arch} is shown in figure \ref{Tented fingerprint}, and figure \ref{Tented NBH fingerprint} illustrates a focused attention around the origin in the phase portrait and the delta in the tented arch fingerprint image. It is clear that the smooth flow and the arrangement of the ridges are similar.
\subsection{Strong arch fingerprint }
\begin{figure}[ht]
\centering
\begin{minipage}[c]{0.35\linewidth}
\centering
\includegraphics[width=2in]{strong.png}
\end{minipage}%
\begin{minipage}[c]{0.35\linewidth}
\centering
\includegraphics[width=1.9in]{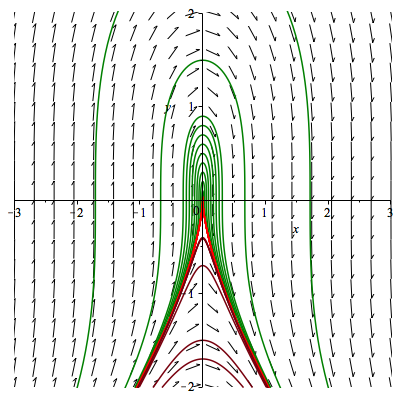}
\end{minipage}%
\caption{Simulation of the strong arch fingerprint}
\label{Strong fingerprint}
\end{figure}

\begin{figure}[ht]
\centering
\begin{minipage}[c]{0.3\linewidth}
\centering
\includegraphics[width= 1.3 in]{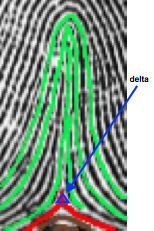}\\ \small (a)
\end{minipage}%
\begin{minipage}[c]{0.4\linewidth}
\centering
\includegraphics[width= 2.3 in]{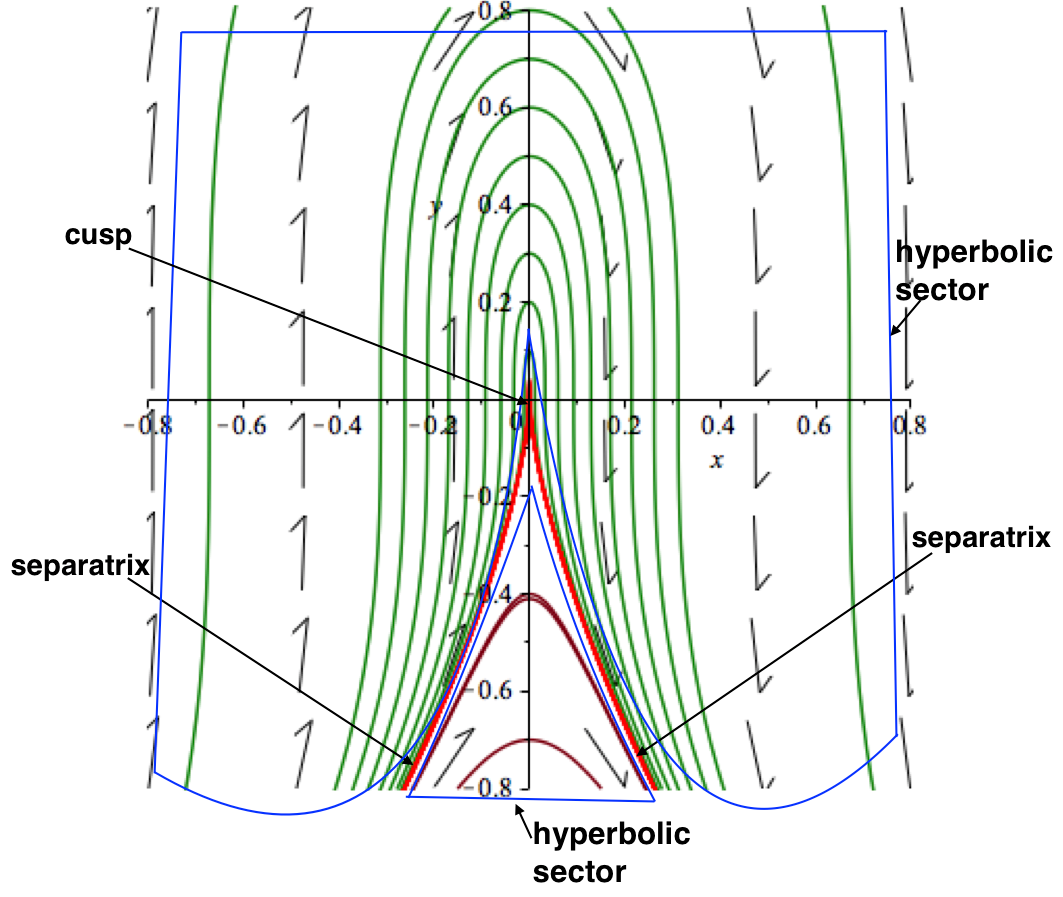}\\ \small (b)
\end{minipage}%
\caption{(a) Neighborhood of the delta in the image of the strong arch fingerprint and (b) Neighborhood of the origin in the phase portrait of example  \ref{example Strong arch}    .}
\label{Strong NBH fingerprint}
\end{figure}

\hspace{\parindent} If $\theta$ grows up enough around  the value $\theta =5$, the flow stretched longer vertically in the middle with a cusp at the origin making acute angle. This flow agrees with the general shape of  the strong arch, for example consider system \eqref{FP-Arch-sys} with  $\theta = 5$ we get the  system \eqref{theta = 5}.
\begin{example} \label{example Strong arch}
\begin{equation} \label{theta = 5}
   \begin{aligned} 
      \frac{dx}{dt} &= y^2,
\\
\frac{dy}{dt}  &= -5x.
\end{aligned} 
\end{equation}
\end{example}
In figure \ref{Strong fingerprint} we see the matching between the image of strong arch and the phase portrait of example \eqref{example Strong arch}, and in figure \ref{Strong NBH fingerprint} we see how the flow is narrowed about the origin point  in a sharp vertical way to be identical with the ridges around the delta in the strong  arch fingerprint image.
\section{Conclusion} \label{FP-Arch conculusion}
\hspace{\parindent}
The general shape of the flow of the considered  dynamical system and  the shape of the ridges in the images of the categories of the arch fingerprint are almost identical to each other.
\bibliographystyle{unsrt}  


\end{document}